# Optimization of Express Train Service Network: Under the Competition of Highway Transportation

Boliang Lin[*]

School of Traffic and Transportation, Beijing Jiaotong University, Beijing 100044, China

**Abstract:** In order to reduce the carbon emission, the related government departments encourage road freights to be transferred more by railway transportation. In China freight transport system, the road transportation is usually responsible for the freights that are in a short distance or the ones with high value-added. To transfer more high value-added freights from highway to railway, except the transportation expenses of railway have an advantage over the road, the transportation time is of certain competitive force as well. Therefore, it is very essential for railway to provide freight train products that are of competitive power. Under such circumstance, a multi-objective programming model of optimizing the rail express train network is devised in this work on the basis of taking both road and railway transportation modes into account. The aims of optimization are to minimize the operation costs of rail trains, and to maximize the railway transport revenue. In a network with a given set of express train services, either the all-or-nothing (AON) method or the logit model can be employed when assigning high value-added freights. These two flow assignment patterns are investigated in this work.

## 1. Introduction

In recent years, with the rapid development of the global economy, the consumption of fossil energy such as coal, oil and natural gas for industrial production and transportation has increased dramatically. At the same time, due for the destruction of the environment, its ability to purify the air is gradually weakened, resulting in an increasing amount of carbon dioxide in the atmosphere. Therefore, enough attention should be paid for to the consequences of global warming and effective measures must be taken to curb it.

Faced with increasing greenhouse gases, reducing carbon emissions and developing green economy have become the consensus of governments around the world. For example, the European Union issued a white paper, titled ''Roadmap to a Single European Transport Area -Towards a competitive and resource efficient transport system,'' stating that 30% of road freight over 300 km should shift to other modes such as rail or waterborne transport by 2030, and more than 50% by 2050 (European Commission [1]). The Ministry of Transport of China also recommends that the railways and waterways should support more freight transportation loads.

---

[*] Corresponding author. Email: bllin@bjtu.edu.cn

In order to improve the railway freight transport, it is necessary to make the most of the comparative advantages of railway. Road freight transport is convenient and efficient. Through continuous reforms, road transportation can meet the requirements of modern logistics development and the transportation needs of high value-added goods. The convenience and efficiency of road transportation resulting in a significant reduction in railway freight market share. Therefore, as one of the largest and busiest railway systems, the China Railway has focused on the development of scheduled train services in past years in an effort to obtain a higher market share in high-value freight transportation.

The scheduled trains in China can be divided into the following three levels according to the operating speed: 160 km/h, 120 km/h and 80 km/h (Lin et al. [2]). The high operating speed of scheduled trains can guarantee the delivery time of the shipments. However, according to the current practice, some train schedule plan is usually made on the basis of historical data and experience. Additionally, some of the train services attract fewer freight flows than expected. In this situation, it is necessary to improve the railway freight transport service to attract more freight flow from the road.

This paper takes highway transportation into consideration, optimizing the flow assignment and concentrating on the block-swap plan of express trains. A bi-level multi-objective programming model is proposed to mathematically describe the problem. The upper-level objective aims to minimize the operation costs of rail trains, and to maximize the railway transport revenue. And the lower-level formulation is constructed to determine the flow assignment and the rail transportation expense, also the train service sequences of high value-added shipments can be settled if the rail transportation mode is selected.

The reminder of this paper is structured as follows: Section 2 presents a comprehensive literature review on related studies. In Section 3, we introduce the optimization of express train network considering the competition from highway transportation. In Section 4, we analyze the optimization objective and constraints, and a multi-objective programming model is proposed. Finally, conclusions are drawn in Section 5.

## 2. Literature Review

To subsequently shift the freights from road to railway, the fundamental issue is to study the choice behavior of transportation modes. Many experts and scholars have carried out researches from this aspect. Tsamboulas et al. [3] proposed a methodology with the necessary tools to assess the potential of a specific policy measure to produce a modal shift in favor of intermodal transport. Crisalli et al. [4] presented a methodology to evaluate rail-road freight policies, which used a specific mode-service choice model to share the freight demand among alternatives and a service network design model to identify new rail-road freight services. Santos et al. [5] discussed the impact of three freight transport policies aiming to promote railroad transport in Europe. These works have provided valuable reference on the measures to shift the freights from road to railway.

Furthermore, the theory of traffic assignment, as the foundation of freight shift, is a classic problem in the field of transportation. Daganzo et al. [6] proposed a quantitative evaluation of probabilistic traffic assignment models and an alternate formulation. Tamin et al. [7] carried out using two assignment techniques, all-or-nothing and the stochastic method in determining the routes taken through the network. The traffic assignment models of comprehensive transportation system was established by Fan et al. [8] using the analytic hierarchy process method and also presented the process of the traffic assignment. The research results of road traffic assignment are difficult to be directly applied to the railway system because they are quite different in terms of transport patterns and transport organization. The exclusive constraints of single flow and tree path limits should be considered for the railway flow assignment problem. Therefore, this kind of model is usually a nonlinear 0-1 programming model. And when the problem scale is large, it is difficult to solve this model. Many researchers thereby adopt similar methods to deal with the flow assignment problem of railway transportation network. Wang et al. [9] investigated the multi-modal express shipment network routing problem based on the concept of express shipment service level, and corresponding model was proposed according to the *k*-shortest path set between each OD. Wang et al. [10] investigated the freight mode choice between truck and rail for the shipments originated from Maryland to other states by using binary probit and logit models. Chen et al. [11, 12] studied the integrated optimization of transfer of freight flow in the land transportation system and the transportation flow assignment, and the corresponding 0-1 mathematic model was put forward. Liu et al. [13] proposed the basic principle of the freight flow assignment, based on which a user equilibrium model for freight flow assignment is formulated. Moreover, sensitivity analysis on the expenditure of road and railway transportation was made to reveal their quantitative influence on the flow distribution, which can provide theoretical method for estimating the measures for diverting bulk freight flow from road to railway.

Optimization of the service network is a powerful way to implement the shift of road flows to railway. Service network design is a complex problem in most cases. It involves many selection (e.g., services, blocks, schedules) and routing (demand loads) decisions. An integrated model is required to address this problem in a comprehensive way. Crainic et al. [14] presented a general modeling framework for the service network design problem for multimode multicommodity freight transportation in the case of a single authority controlling both the service network and the movements of goods. Crainic [15] introduced a new classification of service network design problems and formulations that emphasized the functionality of the formulation rather than the transportation mode to which it is applied. Lulli et al. [16] considered empty cars as a commodity in the mathematical model to design the network of services. Lin et al. [17] modeled the train service problem as a network design problem, and applied a bi-level programming model to the problem,. Zhu et al. [18] addressed the scheduled service network design problem for freight rail transportation, integrating service selection and scheduling, car classification and blocking, train makeup, and routing of time-dependent customer shipments.

There is a rich body of researches that deal with the design of the railway network to promote the efficiency of the transport. The network design problem was mainly aimed at urban traffic at first. For example, Leblanc [19] addressed the problem of determining which links should be improved in an urban road network so that total congestion in the city was minimized. The research on railway network planning started relatively late. Lin et al. [20] investigated the optimization decision for railroad network design problems to determine where a new railroad line to be constructed or an old railroad to be strengthened from view point of project investment and cost of routing cars. They routed cars along the tree-shape for cargo flow assignment, and also included the passenger corridor or high-speed railway project choices in the model. Lin et al. [21] analyzed the carbon emissions of railway and highway transport and estimated the environmental benefit of building a railway, and constructed a bi-level programming model for railway network design. railway network design problem under the background of reducing carbon emissions This work has a strong guiding significance for improving the decision-making of the railway construction in the context of energy conservation and emission reduction.

## 3. Problem Description

Before introducing the flow assignment method, the express train service network should be designed appropriately. The design principle is to maximum the railway profit, which consist of transportation expenses and costs of dispatching express trains. In this study, it is assumed that the service network of highway is preset and the cost of highway service between each node pair is predefined. Besides, the OD matrix of highway and railway demands is given. The service network design problem of the express train (SNDET) and flow assignment are illustrated by the following example.

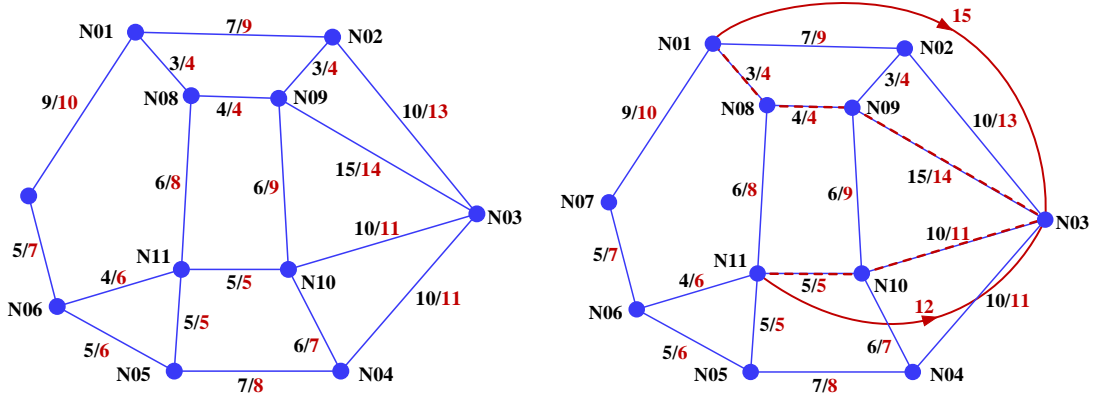

Figure 1. Illustration of service network.   Figure 2. Design of express train service network.

For the highway transportation, a direct service is provided between each node pair. And for the railway transportation, a regular train service is provided only between two adjacent nodes. The limit of delivery time of each demand is of course considered in the high value-added freight transportation, which is the motivation to improve the express train service network. The numbers around arcs in Figure 1a refer to the transportation time of the services of railway and highway. For example, the transportation time of railway service from node N01 to N02 is 9 hours and that of highway service is 7 hours. Besides, for the highway transportation, if the origin and destination of a service are not adjacent, the transportation time is set to the sum of time of intermediate services minus 1 hour multiplied by the number of intermediate nodes. For example, the time of highway service from N07 to N02 is 15 hours.

We assume that there is a demand $F_{73}$ from N07 to N03 with the delivery due time 25 hours. The potential highway transportation plans are as follows:

(i) N07→N01→(N08)→(N09)→N03: 25hr
(ii) N07→(N06)→(N05)→(N04)→N03: 24hr
(iii) N07→(N06)→(N11)→(N10)→N03: 21hr
(iv) N07→N06→N11→N10→N03: 24hr

Nevertheless, there has no qualified railway path for $F_{73}$ based on the current regular train services. Therefore, two express train services are emerged in the railway service network as shown in Figure 2, which are service from N01 to N03 and service from N11 to N03, respectively. The paths of these two services are marked with red dash lines. Under such circumstance, the possible railway transportation plans are given as:

(i) N07→(N01)→(N02)→N03: 25hr
(ii) N07→N06→N11→(N10)→N03: 25hr

The capacities of service arcs of highway network or the regular train network are assumed to be always sufficient for the high value-added freights. The assignment of the flows thereby is only limited by the delivery time. After designing the express train service network, the demands of both railway and highway can be assigned by the all-or-nothing method or the logit model, the formulation of which will be introduced hereinafter.

## 4. Mathematical Model

In this section, we first introduce the notations used in this paper. After that, a bi-level multi-objective programming model is proposed for the service network design problem of the express trains under the competition of highway transportation.

### 4.1 Notations

The notations used in this paper are listed in Table 1.

Table 1 Notations used in this paper.

| Symbol | Definition |
|---|---|
| Sets | |
| $V$ | Set of stations in a rail network; |
| $S^{arc}$ | Set of arcs in a transportation network; |
| $S^{class}$ | Set of express train classes; |
| $S^{flow}$ | Set of origins and destinations of demands; |
| $Ser(i,j,d)$ | Set of all block-swap plans of $d$-class express trains from station $i$ to $j$; |
| $R_{st}$ | Ser of potential paths of demand $F_{st}$ if it is shipped by railway. |
| Parameters | |
| $C_{ij}^d$ | Fixed cost incurred by providing $d$-class express trains from station $i$ to $j$; |
| $T_{ij}^d$ | Period of $d$-class express trains from station $i$ to $j$; (hour) |
| $L_m$ | General distance of the $m$th arc in the transportation service network; (km) |
| $\lambda_d$ | Unit kilometer cost of one $d$-class express train; |
| $\tau_k$ | Block-swap cost at station $k$ of express train; |
| $F_{st}$ | Traffic demand of the O-D pair origins at station $s$ and is destined to station $t$; |
| $T_{st}$ | Delivery due time of demand $F_{st}$; |
| $t_{st}^l$ | Travel time of the $l$th path by railway; |
| $c_{st}^l$ | Fixed cost of one car if demand $F_{st}$ is shipped on the $l$th path by railway; |
| $c_{st}^{highway}$ | Fixed cost of one truck if demand $F_{st}$ is shipped by highway; |
| $b_m$ | Capacity of service arc $m$ in the railway network. |
| Variables | |
| $y_{ijd}^{\mathbb{P}}$ | Decision variables, it takes value one if $d$-class express trains from station $i$ to $j$ are provided with the block-swap plan $\mathbb{P}$; Otherwise, it is zero; |
| $x_{st}^{highway}$ | Decision variables, it takes value one if demand $F_{st}$ is shipped by highway; Otherwise, it is zero; |

| $x_{st}^l$ | Decision variables, it takes value one if demand $F_{st}$ is delivered by the $l$th path of the railway service network; Otherwise, it is zero; |
|---|---|

### 4.2 Model Formulation

The service network design problem of the express trains considering highway transportation can be formulated as a bi-level multi-objective programming model whose objective function and constraints are written as follows:

*Upper-level programming model (UM):*

$$\min Z_1 = \sum_{i \in V} \sum_{j \in V} \sum_{d \in S^{class}} \sum_{\mathbb{P} \in Ser(i,j,d)} \frac{24}{T_{ij}^d} (C_{ij}^d + \sum_{m \in \mathbb{P}} \lambda_d L_m + \sum_{k \in \mathbb{P}} \tau_k) y_{ijd}^{\mathbb{P}} \quad (1)$$

$$\max Z_2 = \sum_{(s,t) \in S^{flow}} C_{st}^{Rail}(Y) f_{st}^{Rail} \quad (2)$$

*Subject to:*

$$y_{ijd}^{\mathbb{P}} \in \{0,1\} \quad \forall i,j \in V, d \in S^{class}, \mathbb{P} \in Ser(i,j,d) \quad (3)$$

The problem involves both the railway and highway transportations, the original demands thereby contains the freights from railway and highway. More importantly, it is assumed that the capacities of service arcs of either the railway network or the highway network are always sufficient for high value-added freights. The objective functions of UM aim at minimize the fixed costs, the general distance costs and block-swap costs of provided express trains, and maximum the railway revenue at the same time. In the equation (2), $C_{st}^{Rail}$ denotes the transportation expenses from station $s$ to $t$, including the operation costs incurred by loading and unloading at stations and the profits of service charges of short-distance logistics. Additionally, $f_{st}^{Rail}$ denotes the demands that are shipped by railway. The values of $C_{st}^{Rail}$ and $f_{st}^{Rail}$ are determined by the lower-level programming model:

*lower-level programming model (LM):*

$$\min \sum_{(s,t) \in S^{flow}} (c_{st}^{highway} x_{st}^{highway} + \sum_{l \in R_{st}} c_{st}^l x_{st}^l) F_{st} \quad (4)$$

*Subject to:*

$$x_{st}^{highway} + \sum_{l \in R_{st}} x_{st}^l = 1 \quad (5)$$

$$\sum_{s \in V} \sum_{t \in V} \sum_{l \in R_{st}} F_{st} x_{st}^l a_{st}^{ml} \leq b_m \quad \forall m \in S^{arc} \quad (6)$$

$$t_{st}^l x_{st}^l \leq T_{st} \quad \forall s,t \in V, l \in R_{st} \quad (7)$$

$$C_{st}^{Rail} = \sum_{l \in R_{st}} c_{st}^l x_{st}^l \quad (8)$$

$$f_{st}^{\text{Rail}} = F_{st} \sum_{l \in R_{st}} c_{st}^l x_{st}^l \tag{9}$$

$$x_{st}^l \in \{0,1\} \quad \forall (s,t) \in S^{\text{flow}}, l \in R_{st} \tag{10}$$

$$x_{st}^{\text{highway}} \in \{0,1\} \quad \forall (s,t) \in S^{\text{flow}} \tag{11}$$

The objective function of LM aims to minimize the costs incurred by transporting demands $F_{st}$, which is assumed to be fixed involving the costs of fuel consumption, crew employment and depreciation, etc. Constraint (5) ensures that only one transportation mode can be used to ship a certain demand. Constraint (6) guarantees the volumes of traffic flows shipped on arc $m$ would not exceed the capacity of it. Since we assume that the capacities of service arcs are always sufficient in the railway network, then $b_m$ can be defined as:

$$b_m = \begin{cases} \dfrac{24}{T_{ij}^d} \Gamma_{ijd}^{\text{TrainSize}} & \text{if } y_{ijd}^{\mathbb{P}} = 1, \ \forall i,j \in V, d \in S^{\text{Degree}}, \ \mathbb{P} \in Ser(i,j,d) \\ M & \text{otherwise} \end{cases} \quad \forall m \in S^{\text{arc}} \tag{12}$$

where $\Gamma_{ijd}^{\text{TrainSize}}$ is the size of the $d$-class express train dispatched from station $i$ to $j$; Constraint (7) ensures the delivery time would not exceed the limit delivery time of demand $F_{st}$ if it is shipped by railway. The variables $C_{st}^{\text{Rail}}$ and $f_{st}^{\text{Rail}}$ are formulated by equations (8) and (9). And equations (10) and (11) are the binary constraints on variables.

**4.3 Flow Assignment**

Two typical flow assignment approaches are analyzed based on the solution $X$ of the SNDET model, including the all-or-nothing (AON) method and the logit method.

*All-or-nothing (AON) method*

The most common approach to traffic assignment is the "all-or-nothing" method. It is assumed that all drivers wish to minimize their perceived travel costs, all the drivers perceive these costs in the same way and that these costs do not depend on flow level (Tamin, Willumsen [7]). It is the fastest and simplest assignment technique. Generally, this assignment is considered to be unrealistic for many networks because it does not consider variations in the effects of congestion. However, in this work, the congestion is not taken into account because the transportation capacities of high value-added freights are always assumed to be sufficient.

To employ the AON method on both the railway and highway network, the general cost should be defined first:

$$c_{st} = \gamma \mathbb{T}_{st} + \Phi_{st} \quad \forall (s,t) \in S^{\text{flow}} \tag{13}$$

where $c_{st}$ is the general cost of trip from origin $s$ to destination $t$ of certain transportation mode, $\mathbb{T}_{st}$ (hour) is the transportation time of this trip, parameter $\gamma$ is used to transfer hours to CNY, and $\Phi_{st}$ (CNY) is the transportation expenses.

The general cost of railway transportation is based on the solution $X$ of the

SNDET model:

$$c_{st}^{Rail}(X) = \gamma \mathbb{T}_{st}^{Rail}(X) + \Phi_{st}^{Rail}(X) \qquad \forall (s,t) \in S^{flow} \tag{14}$$

According to the AON method, when the general cost satisfies the condition $c_{st}^{Rail}(X) < c_{st}^{Highway} - \Delta$, the flow from origin $s$ to destination $t$ should be transported by the railway.

*Logit method*

The model used for the logit method is the discrete choice model which has been by far most widely used in the field of passenger transport choice studies (Wang, Ding, Liu, Xie [10]). The logit model assumes logistic distribution for the error terms. The utility function of the logit model can be given as:

$$U_{st} = V_{st} + \varepsilon_{st} = -(\gamma \mathbb{T}_{st} + \Phi_{st}) + \varepsilon_{st} \qquad \forall (s,t) \in S^{flow} \tag{15}$$

where $U_{st}$ is the utility of the demand $F_{st}$ for alternative transportation mode; $V_{st}$ is the observable portion of the utility; $\varepsilon_{st}$ is the unobservable component of the utility.

The form for the logit choice probability of demand $F_{st}$ selects railway transportation can be expressed as follows:

$$P_{st}^{Rail} = e^{U_{st}^{Rail}(X)} \Big/ \left( e^{U_{st}^{Rail}(X)} + e^{U_{st}^{Highway}} \right) \qquad \forall (s,t) \in S^{flow} \tag{16}$$

## 5. Conclusion

This paper investigated the service network design problem of express trains and flow assignment by the all-or-nothing method or the logit model. A bi-level multi-objective programming model is proposed for the problem. The objectives of the upper programming model aim to minimize the operation costs of rail trains, and to maximize the railway transport revenue, while the objective of the lower programming model is to minimize the costs incurred by transporting demands. The delivery time is considered as a constraint in the lower model. Furthermore, the general cost functions of railway and highway transportations are formulated according to the transportation time and expenses, which are employed in the utility function of the logit model as well. Our future work is to test and apply the proposed model to certain problem instances.